\documentclass[11pt]{article}
\usepackage{mathrsfs, cite}
\usepackage{amsmath,amsfonts,amssymb,rotating,amsthm}
\usepackage{psfrag,eepic,color}
\usepackage{framed}

\def\qed{\nopagebreak\hfill{\rule{4pt}{7pt}}}
\def\proof{\noindent {\it{Proof.} \hskip 2pt}}
\newtheorem{theo}{Theorem}[section]

\newtheorem{coro}[theo]{Corollary}

\newtheorem{ex}{Example}[section]
\theoremstyle{remark}

\begin{document}
\begin{center}
{\Large \bf Gosper Summability of Rational Multiples of \\[5pt]
Hypergeometric Terms
}
\end{center}

\begin{center}
{\large Qing-hu Hou and Guo-jie Li}\\[9pt]
School of Mathematics\\
Tianjin University\\
Tianjin 300072, China\\[6pt]
{\tt
qh\_hou@tju.edu.cn}, {\tt 1017233006@tju.edu.cn}
\end{center}

\vspace{0.3cm} \noindent{\bf Abstract.}
By telescoping method, Sun gave some hypergeometric series whose sums are related to $\pi$ recently. We investigate these series from the point of view of Gosper's algorithm. Given a hypergeometric term $t_k$, we consider the Gosper summability of $r(k)t_k$ for $r(k)$ being a rational function of $k$. We give an upper bound and a lower bound on the degree of the numerator of $r(k)$ such that $r(k)t_k$ is Gosper summable. We also show that the denominator of the $r(k)$ can read from the Gosper representation of $t_{k+1}/t_k$. Based on these results, we give a systematic method to construct series whose sums can be derived from the known ones. We also illustrated the corresponding super-congruences and the $q$-analogue of the approach.

\vskip 15pt


\vskip 15pt


\section{Introduction}
In \cite{Sun11}, Sun derived several identities involving $\pi$ by telescoping method. For example, he deduced
\begin{equation}\label{Sun1.1}
\sum_{k=0}^\infty \frac{k(4k-1){2k \choose k}^3 }{(2k-1)^2 (-64)^k} = - \frac{1}{\pi},
\end{equation}
from Bauer's series
\[
\sum_{k=0}^\infty (4k+1) \frac{{2k \choose k}^3 }{(-64)^k} = \frac{2}{\pi}
\]
and the telescoping sum
\[
\sum_{k=0}^n \frac{(16k^3-4k^2-2k+1){2k \choose k}^2}{(2k-1)^2 (-64)^k} = \frac{8(2n+1)}{(-64)^n}{2n \choose n}^3.
\]
We aim to give a systematic method to construct series like \eqref{Sun1.1}. This motivates us to consider the following problem: Given a hypergeometric term $t_k$, for which rational functions $r(k)$ we have $r(k)t_k$ is Gosper summable?

Recall that a hypergeometric term $t_k$ is Gopser summable \cite{Gosper} (in brief, summable) if there exists a hypergeometric term $z_k$ such that
\[
t_k = z_{k+1} - z_k.
\]
The well-known Gosper's algorithm \cite[Chap.~5]{P-etk} provides us a decision method to find $z_k$ if it exists.
More precisely, suppose that
\[
\frac{t_{k+1}}{t_k} = \frac{a(k)}{b(k)} \frac{c(k+1)}{c(k)},
\]
where $a(k),b(k),c(k)$ are polynomials such that
\[
\gcd(a(k),\, b(k+h)) = 1, \quad \forall\, h \ge 0.
\]
Then $z_k$ exists if and only if there exists a polynomial $x(k)$ such that
\begin{equation}\label{G-eq}
a(k) x(k+1) - b(k-1) x(k) = c(k).
\end{equation}

In many cases, $t_k$ is not summable itself. However, if we multiply it with a suitable rational function $r(k)$, the product $r(k)t_k$ will be summable. By considering the space
\[
\{a(k)x(k+1)-b(k-1)x(k) \,|\, x(k) \in K[k] \},
\]
we derive an upper bound and a lower bound on the degree of the numerator of $r(k)$ in Section~2. By aid of Gosper's algorithm, we also give candidates for the denominator of $r(k)$.

Based on these discussion, we are able to give a systematic method to construct series of the form \eqref{Sun1.1} in Section~3. We not only recover the series given in \cite{Sun11}, but also find some new ones. Moreover, we derive some super-congruences for the partial sums of the series and illustrate the $q$-analogue of this approach by an example.

\section{The numerator and the denominator of $r(k)$} \label{sec-th}
Let $t_k$ be a hypergeometric term over a field $K$ of character zero. That is, $t_{k+1}/t_k$ is a rational function of $k$ over $K$. In this section, we will give some properties on the rational fuunction $r(k)$ such that $r(k)t_k$ is summable.

Suppose
\begin{equation}\label{G-form}
\frac{t_{k+1}}{t_k} = \frac{a(k)}{b(k)} \frac{c(k+1)}{c(k)},
\end{equation}
where $a(k),b(k),c(k) \in K[k] \setminus \{0\}$ and
\begin{equation}\label{G-condi}
\gcd(a(k),\, b(k+h)) = 1, \quad \forall\, h \ge 0,
\end{equation}
we say $(a(k),b(k),c(k))$ is a {\it Gosper representation} of $t_k$.

We first consider the case when $r(k)$ is a polynomial. In this case, we write $p(k)$ instead of $r(k)$. To give an upper bound on the degree of $p(k)$, we need the concept of degeneration introduced in \cite{HMZ}. Let
\begin{equation}\label{def-u}
u(k) = a(k)-b(k-1),
\end{equation}
and
\begin{equation}\label{def-d}
d = \max\{ \deg u(k),\ \deg a(k) - 1 \}.
\end{equation}
The pair $(a(k),b(k))$ is said to be {\it degenerated} if
\[
\deg u(k) = \deg a(k) -1 \quad \mbox{and}  -lc\   u(k)/\ lc \quad a(k) \in \mathbb{N},
\]
where $lc\   f(k)$ denotes the leading coefficient of $f(k)$.

An upper bound on the degree $p(k)$ can be given in terms of $d$.

\begin{theo}\label{th-ub}
Let $t_k$ be a hypergeometric term with a Gosper representation $(a(k),b(k),c(k))$. Then there exists a non-zero polynomial $p(k)$ such that $p(k) t_k$ is summable. Moreover, an upper bound on the degree of $p(k)$ is given by
\begin{equation}\label{deg-b}
B =
\begin{cases}
  d+1, & \mbox{if $(a(k),b(k))$ is degenerated or $\deg u(k) < \deg a(k)-1$}, \\
  d, & \mbox{otherwise}.
\end{cases}
\end{equation}
\end{theo}

\proof
Let
\[
S_{a,b} = \{a(k) x(k+1) - b(k-1) x(k) \,|\, x(k) \in K[k] \}.
\]
By Theorem~2.3 of \cite{HMZ}, the dimension of $K[k]/S_{a,b}$ is bounded by $B$.
Therefore, the $B+1$ vectors
\[
\overline{c(k)},\, \overline{k \cdot c(k)},\, \ldots,\, \overline{k^B \cdot c(k)}
\]
are linearly dependent in $K[k]/S_{a,b}$. Hence
there exists a non-zero polynomial $p(k)$ of degree less than or equal to $B$ such that
$c(k) p(k) \in S_{a,b}$, implying that $p(k) t_k$ is summable. \qed

As a corollary, we have
\begin{coro}\label{def-u/v}
Let $t_k$ be a hypergeometric term and $q(k)$ be a non-zero polynomial. Suppose $(a(k),b(k),c(k))$ is a Gosper representation of $t_k/q(k)$ and $B$ is given by \eqref{deg-b}.
Then there exists a non-zero polynomial $p(k)$ of degree no more than $B$ such that $p(k)t_k/q(k)$ is summable.
\end{coro}

\noindent {\it Remark.} When $a(k),b(k)$ are shift-free, i.e., $\gcd(a(k),\, b(k+h))=1$ for any $h \in \mathbb{Z}$, Chen, Huang, Kauers and Li \cite{Chen} showed that one may take $d$ as the upper bound.

The following examples show that when $a(k),b(k)$ are not shift-free, we need take $d+1$ as the upper bound.
\begin{ex}
Let
\[
t_k = \frac{k}{(k+1)^2(k+2)}.
\]
A Gosper representation of $t_k$ is given by
\[
\frac{a(k)}{b(k)} \cdot \frac{c(k+1)}{c(k)} = \frac{(k+1)^2}{(k+2)(k+3)} \cdot \frac{k+1}{k}.
\]
One sees that $(a(k),b(k))$ is degenerated and $d=1$. But there are no polynomials $p(k)$ of degree less than or equal $1$ such that $p(k)t_k$ is summable.

Let
\[
t_k = \frac{1}{(k+1)^2}.
\]
A Gosper representation of $t_k$ is given by
\[
\frac{a(k)}{b(k)} \cdot \frac{c(k+1)}{c(k)} = \frac{(k+1)^2}{(k+2)^2} \cdot \frac{1}{1}.
\]
One sees that $\deg u(k)<\deg a(k)-1$ and $d=1$. But there are no polynomials $p(k)$ of degree less than or equal $1$ such that $p(k)t_k$ is summable.
\end{ex}

Of course, it is possible that there exists a polynomial $p(k)$ of degree less than $B$ such that $p(k)t_k$ is summable.

\begin{ex}\label{ex-deg}
Let
\[
t_k=\frac{\binom{2k}{k}^4}{(2k-1)^4 256^k}.
\]
A Gosper representation of $t_k$ is given by
\[
\frac{a(k)}{b(k)} \cdot \frac{c(k+1)}{c(k)}=\frac{(2k-1)^4}{16(k+1)^4} \cdot \frac{1}{1}.
\]
We find that $d=3$ and $B=4$, while $(4k-1)t_k$ is summable, in which the degree of $p(k)$ is much smaller than the upper bound.
\end{ex}

In some cases, the upper bound $B$ is also a lower bound.
\begin{theo}\label{th-lb}
Let $t_k$ be a hypergeometric term with a Gosper representation $(a(k),b(k),c(k))$ and $u(k)$ is given by \eqref{def-u}. Suppose that
\begin{equation}\label{cond}
\deg u(k) = \max\{\deg a(k),\, \deg b(k) \} \quad \mbox{and} \quad c(k)=1.
\end{equation}
If $p(k)t_k$ is summable, then
\[
\deg p(k) \ge \max\{\deg a(k),\, \deg b(k) \}.
\]
\end{theo}

\proof Notice that
\[
\deg u(k) = \max\{\deg a(k),\, \deg b(k) \}
\]
implies that for any polynomial $x(k)$,
\[
\deg \big( a(k)x(k+1)-b(k-1)x(k) \big) = \max\{\deg a(k),\, \deg b(k) \} + \deg x(k).
\]
While $p(k)t_k$ is summable if and only if there exists $x(k)$ such that
\[
a(k)x(k+1)-b(k-1)x(k) = p(k).
\]
Therefore,
\[
\deg p(k) \ge \max\{\deg a(k),\, \deg b(k) \} = \deg u(k). \tag*{\qed}
\]

Notice that the condition  $\deg u(k) = \max\{\deg a(k),\, \deg b(k) \}$ implies that $B=d$ in Theorem~\ref{th-ub}.
\begin{ex}
Let
\[
t_k=\frac{\binom{2k}{k}^2}{64^k}.
\]
A Gosper representation of $t_k$ is given by
\[
\frac{a(k)}{b(k)}\cdot \frac{c(k+1)}{c(k)}=\frac{(2k+1)^2}{16(k+1)^2} \cdot \frac{1}{1}.
\]
We find that $B=d=2$ and the quadratic polynomial $p(k)=1+4k-12k^2$ satisfies the condition that  $p(k)t_k$ is summable.
\end{ex}

Example~\ref{ex-deg} shows that the condition $\deg u(k) = \max\{\deg a(k), \deg b(k) \}$ is necessary. The following example indicates that the condition $c(k) \not= 1$ is also necessary.

\begin{ex}
Let
\[
t_k=(3k+2){2k \choose k}.
\]
A Gosper representation of $t_k$ is given by
\[
\frac{a(k)}{b(k)}\cdot \frac{c(k+1)}{c(k)}=\frac{2(2k+1)}{(k+1)} \cdot \frac{3k+5}{3k+2}.
\]
It is easy to check that $B=d=1$, but $t_k$ itself is summable.
\end{ex}

Now we consider the denominator of $r(k)$ such that $r(k)t_k$ is summable.
\begin{theo}\label{th-den}
Let $t_k$ be a hypergeometric term with a Gosper representation $(a(k),b(k),c(k))$ and $p(k),q(k)$ are two polynomials such that $p(k)t_k/q(k)$ is summable. Suppose that
\begin{equation}\label{sfree}
\gcd(q(k),\, a(k-1-h)) = \gcd(q(k),\, b(k+h)) = 1, \quad \forall\, h \in \mathbb{N}
\end{equation}
and
\begin{equation}\label{if-vc}
\gcd(q(k),\, c(k))=1.
\end{equation}
Suppose further that
\begin{equation}\label{if-vv}
\gcd(q(k),\, q(k+1+h)) = 1, \quad \forall\, h \in \mathbb{N}.
\end{equation}
Then we have $q(k) \mid p(k)$.
\end{theo}
\proof
Let $\hat{t}_k = p(k) t_k/q(k)$. We have
\[
\frac{\hat{t}_{k+1}}{\hat{t}_k} = \frac{a(k) q(k)}{b(k) q(k+1)} \cdot \frac{c(k+1)p(k+1)}{c(k)p(k)}.
\]
By \eqref{sfree} and \eqref{if-vv}, we have
\begin{align*}
\gcd(a(k) q(k),\, b(k+h) q(k+1+h)) &=\gcd(a(k),\, b(k+h) q(k+1+h)) \\
&= \gcd(a(k),\, b(k+h)) = 1,  \quad \forall\, h \in \mathbb{N}.
\end{align*}
Therefore,
$(a(k) q(k), b(k) q(k+1), c(k)p(k))$ is a Gosper representation of $p(k)t_k/q(k)$. Hence $\hat{t}_k$ is summable if and only if there exists a polynomial $x(k)$ such that
\[
q(k)\big( a(k)x(k+1)-b(k-1)x(k) \big) = c(k) p(k).
\]
By \eqref{if-vc}, we derive that $q(k) \mid p(k)$. \qed

The following example shows that the condition \eqref{if-vv} is necessary.
\begin{ex}
Let
\[
t_k=\frac{16^k}{\binom{2k}{k}^2}.
\]
A Gosper representation of $t_k$ is given by
\[
\frac{a(k)}{b(k)}\cdot \frac{c(k+1)}{c(k)}=\frac{4(k+1)^2}{(2k+1)^2}\cdot \frac{1}{1}.
\]
Let
\[
r(k)=\frac{p(k)}{q(k)}=\frac{1-40k-32k^2}{(4k+1)(4k+5)}.
\]
It is easy to check that $r(k)t_k$ is summable and $q(k)$ satisfies \eqref{sfree} and \eqref{if-vc}.
\end{ex}

Let $r(k)$ be a rational function such that $r(k)t_k$ is summable. Write $r(k)=p(k)/q(k)$ in reduced form, i.e., $\gcd(p(k),q(k))=1$ and assume that $q(k)$ satisfies \eqref{if-vv}. Then Theorem~\ref{th-den} tells us that $q(k)$ contains a factor of $c(k)$, $a(k-1-h)$, or $b(k+h)$ for some $h \in \mathbb{N}$.

\section{Constructing new hypergeometric series}
In this section, we will use the results in Section~\ref{sec-th} to give a systematic method of constructing series similar to \eqref{Sun1.1}.

Given a hypergeometric series whose sum is known
\[
\sum_{k=0}^\infty t_k = C.
\]

We first compute a Gosper representation $(a(k),b(k),c(k))$ of $t_k$ and make $c(k)$ to be $1$ by setting $\hat{t}_k = t_k/c(k)$. Then we try to find rational functions $r(k)=p(k)/q(k)$ such that $(c(k)+r(k))\hat{t}_k$ is summable. Theorem~\ref{th-den} tells us that a good choice for $q(k)$ is a factor of $a(k-1)b(k)$. Given $q(k)$, we can search for $p(k)$ by applying the Extended Zeilberger's algorithm \cite{CHM} to
\[
t_k,\ \hat{t}_k/q(k),\ k \hat{t}_k/q(k),\ \ldots,\ k^m \hat{t}_k/q(k).
\]
In fact, the algorithm finds $a_0,a_1,\ldots,a_m$ and $g(k)$ such that
\[
t_k + (a_0+a_1 k+ \cdots a_m k^m) \hat{t}_k/q(k) = g(k+1)-g(k).
\]
Notice that Corollay~\ref{def-u/v} ensures the existence of non-trivial $a_i$'s and Theorem~\ref{th-lb} gives a lower bound for the degree $m$.
Let
\[
p(k) = a_0+a_1 k+ \cdots a_m k^m.
\]
We finally derive a new hypergeometric series
\[
\sum_{k=0}^\infty \frac{p(k)}{q(k)} \frac{t_k}{c(k)} = -C + \lim_{k \to \infty} g(k) - g(0).
\]

Note that when \{$\deg p(k)-\deg q(k)\}<\deg c(k)$, the new series convergent faster than the original one.

We will show the method by the following example.
\begin{ex}\label{ex3.1}
Consider Bauer's series
\[
\sum_{k=0}^{\infty}\frac{(4k+1)\binom{2k}{k}^3}{(-64)^k}=\frac{2}{\pi}.
\]

Let
\[
t_k=\frac{(4k+1)\binom{2k}{k}^3}{(-64)^k}.
\]
A Gosper representation of $t_k$ is given by
\[
\frac{a(k)}{b(k)}\cdot \frac{c(k+1)}{c(k)}=\frac{-(2k+1)^3}{8(k+1)^3}\cdot \frac{4k+5}{4k+1}.
\]
So we set
\[
\hat{t}_k =  \frac{t_k}{4k+1} = \frac{\binom{2k}{k}^3}{(-64)^k}.
\]

We first take $q(k)$ to be the factor $(2k-1)^2$ of $a(k-1)$. In this case, we find that
\[
t_k +\frac{(8k^2-2k)}{(2k-1)^2} \hat{t}_k = g(k+1)-g(k),
\]
with
\begin{equation}\label{g1}
g(k) = -\frac{8 k^3}{(2k-1)^2} \frac{\binom{2k}{k}^3}{(-64)^k}.
\end{equation}
Noting that $\lim_{k \to \infty}g(k)=g(0)=0$, we thus derive
\[
\sum_{k=0}^\infty \frac{(4k-1)k}{(2k-1)^2} \frac{\binom{2k}{k}^3}{(-64)^k} = -\frac{1}{\pi},
\]
which is exactly \eqref{Sun1.1}.

Then we take $q(k) = (k+1)^2$ which is a factor $b(k)$. This time we find that
\[
t_k - \frac{1}{4} \frac{(4k+3)(2k+1)}{(k+1)^2} \hat{t}_k = g(k+1) - g(k),
\]
with
\begin{equation}\label{g2}
g(k) = -2k \frac{\binom{2k}{k}^3}{(-64)^k}.
\end{equation}
We thus derive
\begin{equation}
\sum_{k=0}^{\infty}\frac{(4k+3)(2k+1)}{(k+1)^2}\frac{\binom{2k}{k}^3}{(-64)^k}=\frac{8}{\pi}.
\end{equation}

Finally we take $q(k)=(k+1)(2k-1)$. In this case, we discover the identity
\begin{equation}\label{q=ab}
\sum_{k=0}^{\infty}\frac{4k+1}{(k+1)(2k-1)}\frac{\binom{2k}{k}^3}{(-64)^k}=-\frac{4}{\pi}.
\end{equation}

Similarly, when $q(k)$ is of degree $3$, we get two more identities
\begin{align*}
& \sum_{k=0}^\infty \frac{4k-1}{(2k-1)^3} \frac{\binom{2k}{k}^3}{(-64)^k} = \frac{2}{\pi}, \\[5pt]
& \sum_{k=0}^\infty \frac{4k+3}{(k+1)^3} \frac{\binom{2k}{k}^3}{(-64)^k} = 8 - \frac{16}{\pi},
\end{align*}
where the constant $8$ comes from $g(0)$.
\end{ex}

Here is another example, in which the new series may not start from the first term of the original one.
\begin{ex}
In 1993, Zeilberger\cite{Zeilber} used the WZ method to show that
\[
\sum_{k=1}^{\infty}\frac{21k-8}{k^3\binom{2k}{k}^3}=\frac{\pi^2}{6}.
\]

In this case, we have
\[
\frac{a(k)}{b(k)}\cdot \frac{c(k+1)}{c(k)}=\frac{k^3}{8(2k+1)^3}\cdot \frac{21k+13}{21k-8}
\]
and
\[
\hat{t}_k=\frac{1}{k^3\binom{2k}{k}^3}.
\]

Taking $q(k)=(2k+1)^2$, we derive
\begin{equation}
\sum_{k=1}^{\infty}\frac{28k^2+31k+8}{(2k+1)^2k^3\binom{2k}{k}^3}=\frac{\pi^2-8}{2}.
\end{equation}
If we take $q(k)=(k-1)^2$, the series should start from $k=2$. We thus derive
\begin{equation}
\sum_{k=2}^{\infty}\frac{7k^2-8k+2}{(k-1)^2 k^3\binom{2k}{k}^3}=\frac{10-\pi^2}{16}.
\end{equation}
\end{ex}

By considering the congruent properties of the partial sums of the series and the boundary values, we are able to derive some super-congruences.

\begin{ex}
The Bauer's series have a nice $p$-adic analogue \cite{Mort}
\[
\sum_{k=0}^{\frac{p-1}{2}} \frac{(4k+1) \binom{2k}{k}^3 }{(-64)^k}\equiv p\left(\frac{-1}{p}\right) \pmod{p^3}.
\]
By Example~\ref{ex3.1}, we have
\[
\frac{(4k+1) \binom{2k}{k}^3 }{(-64)^k} + \frac{(8k^2-2k)}{(2k-1)^2}\frac{\binom{2k}{k}^3 }{(-64)^k}
= g(k+1)-g(k),
\]
where
\[
g(k) = -\frac{8 k^3}{(2k-1)^2} \frac{\binom{2k}{k}^3}{(-64)^k}.
\]
It is known \cite{Morl} that  for prime $p>3$,
\[
\binom{p-1}{\frac{p-1}{2}} \equiv (-1)^{\frac{p-1}{2}} 2^{2p-2} \pmod{p^3}.
\]
Therefore,
\[
g \left(\frac{p+1}{2} \right) = (-64)^{-(p-1)/2} p \binom{p-1}{\frac{p-1}{2}}^3
\equiv p 8^{p-1} \pmod {p^4}.
\]
Summing over $0$ to $(p-1)/2$, we derive that
\begin{equation}
\sum_{k=0}^{\frac{p-1}{2}} \frac{(8k^2-2k)}{(2k-1)^2}\frac{\binom{2k}{k}^3 }{(-64)^k}
= - p\left(\frac{-1}{p}\right) + p 8^{p-1} \pmod{p^3}.
\end{equation}
Note that when $p=3$, the above congruence also holds.

In a similar way, we derive that
\begin{equation}
\sum_{k=0}^{\frac{p-1}{2}} \frac{(4k+3)(2k+1)}{4 (k+1)^2} \frac{\binom{2k}{k}^3}{(-64)^k}
\equiv p\left(\frac{-1}{p}\right) \pmod{p^3},
\end{equation}
and
\begin{equation}
\sum_{k=0}^{\frac{p-1}{2}} \frac{4k+1}{2(k+1)(2k-1)} \frac{\binom{2k}{k}^3}{(-64)^k}
\equiv p^2 - p \left(\frac{-1}{p}\right) \pmod{p^3}.
\end{equation}
\end{ex}

We remark that all the discussion works smoothly for the $q$-analogue of hypergeometric series. We will give an example on the $q$-case to conclude the paper.
\begin{ex}
A $q$-analogue of Bauer's identity is given by \cite[Eq.(1.5)]{Guojunyi}
\[
\sum_{k=0}^{\infty}\frac{(-1)^k(1-q^{4k+1})q^{k^2}}{1-q}\frac{(q;q^2)_{k}^{3}}{(q^2;q^2)_{k}^{3}}=\frac{(q;q^2)_\infty (q^3;q^2)_\infty}{(q^2;q^2)_{\infty}^{2}}.
\]

Let
\[
t_k=\frac{(-1)^k(1-q^{4k+1})q^{k^2}}{1-q}\frac{(q;q^2)_{k}^{3}}{(q^2;q^2)_{k}^{3}}.
\]
A $q$-Gosper representation of $t_k$ is given by
\[
\frac{a(k)}{b(k)}\cdot \frac{c(k+1)}{c(k)} = \frac{-q^{2k+1}(1-q^{2k+1})^3}{(1-q^{2k+2})^3}\cdot \frac{1-q^{4k+5}}{1-q^{4k+1}},
\]
so that
\[
\hat{t}_k = \frac{(-1)^kq^{k^2}}{1-q} \frac{(q;q^2)_{k}^{3}}{(q^2;q^2)_{k}^{3}}.
\]

Taking $q(k)=(1-q^{2k-1})^2$, we derive
\begin{equation}\label{q-logue1}
\sum_{k=0}^{\infty} \frac{(-1)^{k}(1-q^{2k})(1-q^{4k-1}) q^{k^2-1}} {(1-q^{2k-1})^2}\frac{(q;q^2)_{k}^{3}}{(q^2;q^2)_{k}^{3}} = -
\frac{(q;q^2)_\infty (q^3;q^2)_\infty}{(q^2;q^2)_{\infty}^{2}},
\end{equation}
which is a $q$-analogue of \eqref{Sun1.1} and Equation~(1.9) of \cite{HS}.

Taking $q(k) = (1-q^{2k-1})(1-q^{2k+2})$, we derive
\begin{equation}\label{q-logue2}
\sum_{k=0}^{\infty}\frac{(-1)^{k} q^{k^2+2k-1} (1-q^{4k+1} )}{(1-q^{2k-1}) (1-q^{2k+2})}\frac{(q;q^2)_{k}^{3}}{(q^2;q^2)_{k}^{3}}= -
\frac{(q^3;q^2)^2_\infty}{(q^2;q^2)_{\infty}^{2}},
\end{equation}
which is a $q$-analogue of \eqref{q=ab}.
\end{ex}

\vskip 0.2cm
\noindent{\bf Acknowledgments.}
This work is supported by the National Natural Science Foundation of China (grants 11771330 and 11921001).

\end{document}